\documentclass[10pt,a4paper]{article}
\usepackage[utf8]{inputenc}
\usepackage[T1]{fontenc}
\usepackage{amsmath}
\usepackage{amssymb}
\usepackage{graphicx}

\newcommand{\ZZ}{{\mathbb Z}}

\newtheorem{cor}{\noindent\rm\bf Corollary}[section]
\newtheorem{thm}{\noindent\rm\bf Theorem}[section]

\newtheorem{lem}{\noindent\rm\bf Lemma}[section]

\begin{document}
	
\title{Value sets of non-permutation polynomials over the residue class rings of integers}
\author{Shikui Shang}
\date{}
\maketitle


\vspace{2mm}

\abstract{In this paper, we study the value sets of non-permutation polynomial functions over the residue class ring $\ZZ/m\ZZ$. When $m=p^r$ is a power of some prime $p$, an upper bound is given for the size of the value set of a polynomial function which is not a permutation. We also show that this upper bound can be achieved by some integral polynomials. Finally, we generalize the results for any positive integer $m$ 
with known prime decomposition.
}



\section{Introduction}	

Assume that $R$ is a unital commutative ring. Let $$f=\sum_{i=0}^da_iX^i\in R[X], \ a_i\in R$$ be a polynomial over $R$ with $\deg f=d$. $f$ induces a function
$\Phi_f:R\rightarrow R, a\mapsto f(a)$. In general, $\Phi_f$ is neither injective nor surjective. Set $V(f)=\text{Im}(\Phi_f)$ is the value set of $f$ on $R$. Then, $\Phi_f$ is surjective if and only if $V(f)=R$.

If $\Phi_f$ is a permutation over $R$, i.e., $\Phi_f$ is both injective and surjective, $f$ is a called a permutation polynomial. We see that $V(f)=R$ is a necessary condition for that $f$ is a permutation polynomial. When $R$ is finite, it is also sufficient.

The most famous finite commutative rings are finite fields, of which the cardinalities are powers of primes. The value set problems for polynomial functions over finite fields are classic topics studied by many papers(e.g.  \cite{CLMS}, \cite{CG-CM}, \cite{CHW}, \cite{DM},\cite{G-CM} and \cite{M}). Besides them, there are several results on the upper bounds of $|V(f)|$ for a non-permutation polynomial $f$ obtained by Wan and his coauthors in \cite{GW}, \cite{MWW}, \cite{W} and \cite{WSC}. All
of these results relate $|V(f)|$ to the degree $d$ of $f$. For example, Wan proved that 
$$|V(f)|\leq q-\lceil\frac{q-1}{d}\rceil.$$
in \cite{W}.

In contrast to finite fields, the study on the value set of polynomial functions over other finite commutative rings are few. In \cite{B}, the value sets of Dickson Polynomials over Galois rings were studied.

The residue class rings of integers are also elementary. Some topics on polynomial functions over residue class rings were studied in \cite{C} and \cite{K}. Specifically, the papers such as \cite{KO}, \cite{MS} and \cite{S} focus on counting the number of polynomial functions.

In this paper, we study the upper bound of the sizes of value sets for non-permutation polynomials over the residue class rings $\ZZ/m\ZZ$. For convenience, we always consider integral polynomials in $\ZZ[X]$. Take a polynomial $f\in\ZZ[X]$. For any $m\geq 1$, the polynomial $(f\bmod m)\in(\ZZ/m\ZZ)[X]$ induces a function on $\ZZ/m\ZZ$. Set
$$N(f,m)=\# V(f\bmod m)$$
is the size of the value set of $(f\bmod\ m)$. Then,
$(f\bmod\ m)$ is a permutation if and only if $N(f,m)=m$.

If $f(\bmod\ m)$ is not a permutation, then $N(f,m)<m$. Define
$$M(m)=\max\{N(f,m)\ |\ (f\bmod\ m) \text{ is not a permutation.}\}$$
is the upper bound of the sizes of value sets for non-permutation polynomials over $\ZZ/m\ZZ$.

When $m=p^r$ is a power of a prime $p$, we can use the $p$-adic lifting method to compute the precise value of $M(p^r)$. That is the main result of this paper.
\begin {thm}\label{thm:1.1}For a prime $p$ and $r\geq 1$, $$M(p^r)=\begin{cases}p-1, &r=1\\
	p^{r-2}(p^2-p+1), &r>1\end{cases}.$$
\end{thm}	

This result tell us that even without the restriction on the degree of polynomials, the value sets of polynomial functions over $\ZZ/p^r\ZZ$ are far from being arbitrary. When $r\geq 2$, the gap of the sizes of the value sets between permutation polynomials and non-permutation polynomials is $p^{r-1}-p^{r-2}$ at least. 

Furthermore, we give the value of $M(m)$ for a general positive integer $m$ with known prime decomposition using Chinese remainder theorem in the last section.

\section{Preliminaries}	

Consider an integral polynomial $f(X)\in\ZZ[X]$ and an integer $m\geq 1$. An important property of $f$ is that for any $a,b\in\ZZ$, $a\equiv b(\bmod\ m)$ implies $f(a)\equiv f(b)(\bmod\ m)$.

Let $f(X)=\sum_{i=0}^da_iX^i\in\ZZ[X]$ of degree $d$.
Set
$$f^{(1)}(X)=f'(X)=\sum_{i=1}^dia_iX^{i-1}$$
is the (formal) derivation of $f(X)$.
Inductively, let
$$f^{(n+1)}(X)=(f^{(n)}(X))'.$$

There exists the expansion formula of Taylor
\begin{equation}f(X+Y)=\sum_{i=0}^d\frac{f^{(i)}(X)}{i!}Y^i.\tag{2.1}\end{equation}
Notice that for any $i\geq0, a\in\ZZ$, $\frac{f^{(i)}(a)}{i!}$ is always an integer.

Applying Hensel's Lemma in $p$-adic analysis, we have  the following well-known theorem(Theorem 123 in \cite{HW}) 
\begin {thm}\label{thm:2.1}(The criterion for permutation polynomial functions)For a prime $p$, an integer $r\geq 2$ and a polynomial $f\in\ZZ[X]$, $f$ is a permutation over $\ZZ/p^r\ZZ$ if and only if the following tow conditions hold:

(1)  $(f\bmod\ p)$ is a a permutation polynomial over $\ZZ/p\ZZ$;

(2)  $f'(a)\not\equiv0(\bmod\ p)$ for any $a\in\ZZ$. $\Box$
\end{thm}	

\begin {cor}\label{cor:2.1}When $r\geq 2$, $f$ is a permutation over $\ZZ/p^r\ZZ$ if and only if $f$ is a permutation over $\ZZ/p^2\ZZ$. $\Box$
\end{cor}

{\bf Remark. } In \cite{R}, R. Rivest also gave a more concrete criterion for permutation polynomial functions when $p=2$.
A polynomial $f(X)=\sum_{i=0}^da_iX^i\in\ZZ[X]$ is a permutation polynomial over $\ZZ/2^r\ZZ$ for $r\geq 2$ if and only if the following
three linear identities over $\ZZ/2\ZZ$ hold
\begin{align*}&a_1\equiv1(\bmod\ 2),\\
	&a_2+a_4+a_6+\cdots\equiv0(\bmod\ 2),\\
	&a_3+a_5+a_7+\cdots\equiv0(\bmod\ 2).
\end{align*}

\

For any $a\in\ZZ\setminus\{0\}$, define $\text{ord}_p(a)=\max\{i\geq 0\ |\ p^i\mid a\}$.
For convenience, we take $\text{ord}_p(0)=\infty$.

The following Lemma is also useful.
\begin {lem}\label{lem:2.1}(A variant of generalized Hensel's Lemma)Assume that $f\in\ZZ[X]$,
$a,b\in\ZZ$, $s=\text{ord}_p(f'(a))$.

(1) If $r\geq 2s$, then $a\equiv b(\bmod\ p^{r-s})$ implies $f(a)\equiv f(b)(\bmod\ p^r)$.

(2) If there exists some $r_0\geq 2s+1$ such that 
$$a\equiv b(\bmod\ p^{r_0-s})\Longleftrightarrow f(a)\equiv f(b)(\bmod\ p^{r_0}),$$ then for any $r\geq r_0$, \begin{equation}a\equiv b(\bmod\ p^{r-s})\Longleftrightarrow f(a)\equiv f(b)(\bmod\ p^r).\tag{2.2}\end{equation}

(3) If $s=0$ and $a\equiv b(\bmod\ p)$, then for all $r\geq 1$,
\begin{equation}a\equiv b(\bmod\ p^r)\Longleftrightarrow f(a)\equiv f(b)(\bmod\ p^r).\tag{2.3}\end{equation} 
\end{lem}	

{\bf Proof. }Since $s=\text{ord}_p(f'(a))$, $f'(a)\equiv0(\bmod\ p^s)$ but $f'(a)\not\equiv0(\bmod\ p^{s+1})$.

(1) Assume that $r\geq 2s$. If $a\equiv b(\bmod\ p^{r-s})$,
there exists some $l\in\ZZ$ such that $b=a+lp^{r-s}$.
When $i\geq 2$,
 $$i(r-s)\geq 2(r-s)\geq r.$$
 
Taking $X=a,Y=lp^{r-s}$ and modulo $p^r$ , we only need the first two terms in $(2.1)$. That is,
$$f(a)=f(b+lp^{r-s})\equiv f(a)+f'(a)lp^{r-s}\equiv f(a)(\bmod\ p^r).$$

(2) We only need to show that $f(a)\equiv f(b)(\bmod\ p^r)$ implies $a\equiv b(\bmod\ p^{r-s})$ for all $r\geq r_0$.

we induct on $r$. For $r=r_0$, it is from the assumption.

If it is true that $$f(a)\equiv f(b)(\bmod\ p^r)\Longrightarrow a\equiv b(\bmod\ p^{r-s}),$$ then for $r+1$, first one has $$f(a)\equiv f(b)(\bmod\ p^{r+1})\Rightarrow f(a)\equiv f(b)(\bmod\ p^r)\Rightarrow a\equiv b(\bmod\ p^{r-s}).$$
 
Set $b=a+lp^{r-s}$ for some $l\in\ZZ$. Since $r\geq r_0\geq 2s+1$, $i(r-s)\geq 2(r-s)\geq r+1$ for $i
\geq 2$.
Modulo $p^{r+1}$, we also have  
$$f(b)\equiv f(a)+f'(a)lp^{r-s}(\bmod\ p^{r+1}).$$
Hence, $f(a)\equiv f(b)(\bmod\ p^{r+1})$ yields
$f'(a)l\equiv0(\bmod\ p^{r+1})$. But $f'(a)\not\equiv 0(\bmod\ p^{s+1})$, we obtain that $l\equiv0(\bmod\ p)$. And, $$b=a+lp^{r-s}\equiv a(\bmod\ p^{r+1-s}).$$

Therefore, (2.2) holds for any $r\geq r_0$.

(3) When $a\equiv b(\bmod\ p)$,
$$a\equiv b(\bmod\ p)\Longleftrightarrow f(a)\equiv f(b)(\bmod\ p)$$
holds unconditionally. Taking $s=0$ and $r_0=1$ in (2), we obtain (3). $\Box$

{\bf Remark.} (1) The assumption on the existence of $r_0$ in (2) is necessary. Otherwise, $f(a)\equiv f(b)(\bmod\ p^r)$ may not imply $a\equiv b(\bmod\ p^{r-s})$ even when $r\geq 2s+1$ holds. A counterexample is as follows. Set $f(X)=X^2+p(p-1)X$, $a=0$. Then, $f'(0)=p(p-1)$ and $s=\text{ord}_p(f'(0))=1$. When $r=3$,
we have that $f(p)=p^3\equiv0\equiv f(0)(\bmod\ p^3)$ but $p\not\equiv 0(\bmod\ p^2)$. Here, $r_0$ can be seen as the base to begin $p$-adic lifting.

(2) Comparing with the original Hensel's lemma, the generalized Hesnel's lemma works for the polynomial functions of which the derivations have higher $p$-orders somewhere. A simple version of the generalized Hensel's lemma can be found in Exercise 6 of Section 1.5 in \cite{Ko}.

\section{The proof of the main result}

In this section, we prove the main theorem of this paper. First, the upper bound of $M(p^r)$ is easily to be obtained by the following two lemmas.

\begin {lem}\label{lem:3.1}Assume that $f\in\ZZ[X]$.
If $f$ is not a permutation over $\ZZ/p\ZZ$, then for any $r\geq 1$
$$N(f,p^r)\leq p^{r-1}(p-1).$$
\end{lem}	
{\bf Proof. }Since $f$ is not a permutation over $\ZZ/p\ZZ$, there is $b_0\in\ZZ/p\ZZ$ such that
$f(a)\not\equiv b_0(\bmod\ p)$ for any $a\in\ZZ$.

For a fixed $r\geq 1$, set
$$T=\{b\in\ZZ/p^r\ZZ\ |\ b\equiv b_0(\bmod\ p)\}.$$
Then, $|T|=p^{r-1}$ and $V(f\bmod\ p^r)\cap T=\phi$. Hence,
$$V(f\bmod\ p^r)\subseteq(\ZZ/p^r\ZZ)\setminus T$$
and $N(f,p^r)\leq p^r-p^{r-1}=p^{r-1}(p-1)$.  $\Box$	

\begin {lem}\label{lem:3.2}Assume $f\in\ZZ[X]$.
If there exists an integer $a\in\ZZ$ such that $f'(a)\equiv0(\bmod\ p)$, then for $r\geq 2$,
$$N(f,p^r)\leq p^{r-2}(p^2-p+1).$$
\end{lem}	
{\bf Proof.} Without lose of generality, we can take an integer $a_0$ in $\{0,1,\cdots,p-1\}$ such that $f'(a_0)\equiv0(\bmod\ p)$. Then, for $a,b\in\ZZ$, if $a\equiv a_0(\bmod\ p)$ and $a\equiv b(\bmod\ p^{r-1})$,
i.e., for some $l\in\ZZ$, $b=a+lp^{r-1}$, we have that
$f'(a)\equiv f'(a_0)(\bmod\ p)$ and 
$$f(b)\equiv f(a)+f'(a)p^{r-1}\equiv f(a)+f'(a_0)lp^{r-1}\equiv f(a)(\bmod\ p^r)$$ since $r\geq 2$.	

Consider the set $R_r=\{0,1,\cdots,p^r-1\}$ of all representives in $\ZZ/p^r\ZZ$. Definte the subset $S_r(a_0)$ of $R_r$ by
$$S_r(a_0)=\{a_0+lp+kp^{r-1}\ |\ 0\leq l\leq p^{r-2}, 1\leq k\leq p-1\}.$$
We have seen that $f(a_0+lp+kp^{r-1})\equiv f(a_0+lp)(\bmod\ p^r)$ for $1\leq k\leq p-1$. But notice that $a_0+lp\notin S_r(a_0)$. One has that
$$f(S_r(a_0))\subseteq f(R_r\setminus S_r(a_0)).$$
Hence,
$$V(f\bmod\ p^r)=f(R_r\setminus S_r(a_0))\cup f(S_r(a_0))=f(R_r\setminus S_r(a_0)).$$
Consequently,
\begin{align*}&N(f,p^r)=|V(f\bmod\ p^r)|=|f(R_r\setminus S_r(a_0))|\\
\leq&|R_r\setminus S_r(a_0))|=|R_r|-|S_r(a_0)|=p^r-p^{r-2}(p-1)\\
=&p^{r-2}(p^2-p+1).\ \ \ \Box \end{align*}

Next, we will show the upper bound of $M(p^r)$ can be achieved by some polynomial function over $\ZZ/p^r\ZZ$. 

\begin{lem}\label{lem:3.3}For $f(X)\in\ZZ[X]$, assume
that $f'(a)\not\equiv0(\bmod\ p^2)$ for any $a\in\ZZ$ and $N(f,p^3)=p^3-p^2+p$, then for any $r\geq 3$,
$$N(f,p^r)=p^{r-2}(p^2-p+1).$$
\end{lem}	
{\bf Proof.} First, by Lemma \ref{lem:3.1}, $N(f,p^3)=p^3-p^2+p$ yields that $f$ is a permutation over $\ZZ/p\ZZ$. As we have seen in the proof of Lemma \ref{lem:3.2}, $N(f,p^3)=p^3-p^2+p$ also means that there exits an $a_0\in\{0,\cdots,p-1\}$ such that $a\equiv a_0(\bmod\ p)\Rightarrow f'(a)\equiv0(\bmod\ p)$ and $a\not\equiv a_0(\bmod\ p)\Rightarrow f'(a)\not\equiv0(\bmod\ p)$ for $a\in\ZZ$. Hence, associated with the assumption of $f'(a)$,
$$\text{ord}_p(f'(a))=\begin{cases}1, &a\equiv a_0(\bmod\ p)\\0, &a\not\equiv a_0(\bmod\ p)\end{cases}.$$

We use the same notions in the proof of Lemma \ref{lem:3.2}, one has that
$$|f(R_3\setminus S_3(a_0))|=N(f,p^3)=p^3-p^2+p=|R_3\setminus S_3(a_0))|,$$
which show that $(f\bmod\ p^3)$ is injective restricted on the set $R_3\setminus S_3(a_0)$.
Therefore, $f(a)\equiv f(b)(\bmod\ p^3)$ holds if and only if one of the following two conditions is satisfied:

(1) $a\equiv a_0(\bmod\ p), \text{ord}_p(f'(a))=1$ and $a\equiv b(\bmod\ p^2)$;

(2) $a\not\equiv a_0(\bmod\ p), \text{ord}_p(f'(a))=0$ and $a\equiv b(\bmod\ p^3)$.

When $a\equiv a_0(\bmod\ p)$, the assumption of Lemma \ref{lem:2.1} (2) is satisfied for $s=1$. And, when $a\not\equiv a_0(\bmod\ p)$, the assumption of Lemma \ref{lem:2.1} (3) is satisfied. We obtain that for $r\geq 3$, $f(a)\equiv f(b)(\bmod\ p^r)$ holds if and only if one of the following two conditions is satisfied:

(1) $a\equiv a_0(\bmod\ p)$ and $a\equiv b(\bmod\ p^{r-1})$;

(2) $a\not\equiv a_0(\bmod\ p)$ and $a\equiv b(\bmod\ p^r)$.

Thus, $(f\bmod\ p^r)$ is injective on $R_r\setminus S_r(a_0)$. By Lemma \ref{lem:3.2},
\begin{align*}&p^{r-2}(p^2-p+1)=|R_r\setminus S_r(a_0)|=|f(R_r\setminus S_r(a_0))|\\
	\leq&|f(R_r)|=N(f,p^r)\leq p^{r-2}(p^2-p+1),\end{align*}
which implies $N(f,p^r)=p^{r-2}(p^2-p+1)$.    $\Box$

\begin {lem}\label{lem:3.4}For a prime $p$, set $f(X)=X^{2p-1}+pX\in\ZZ[X]$. Then, for any $r\geq 2$,
$$N(f,p^r)=p^{r-2}(p^2-p+1).$$
\end{lem}	
{\bf Proof.} First, by Fermat's little theorem, for any $a\in\ZZ$, $a^p\equiv a(\bmod\ p)$ and
$$f(a)=a^{2p-1}-pa\equiv a^{2p-1}\equiv a(\bmod\ p).$$
Meanwhile, $f'(X)=(2p-1)X^{2p-2}+p$,
$$f'(a)=(2p-1)a^{2p-2}+p\equiv -a^{2p-2}\equiv -a^{p-1}(\bmod\ p).$$

For $a\not\equiv 0(\bmod\ p)$, $f'(a)\equiv -1(\bmod\ p)$. And, for $a\equiv 0(\bmod\ p)$, $f'(a)=(2p-1)a^{2p-2}-p\equiv-p\not\equiv 0(\bmod\ p^2)$ since $2p-2\geq2$ for any prime $p$.

We always have that $f'(a)\not\equiv0(\bmod\ p^2)$. By Lemma \ref{lem:3.3}, we only need to show that $N(f,p^r)=p^{r-2}(p^2-p+1)$ for $r=2,3$.

When $r=2$, since $f(a)\equiv a(\bmod\ p)$, $(f\bmod\ p)$ is a permutation over $\ZZ/p\ZZ$. Moreover, the fact that
$$\text{ord}_p(f'(a))=\begin{cases}1, &a\equiv 0(\bmod\ p)\\0, &a\not\equiv 0(\bmod\ p)\end{cases}$$
implies $(f\bmod\ p^2)$ is injective on $R_2\setminus S_2(0)$, where $S_2(0)=\{lp\ |\ 1\leq l\leq p-1\}$.
Hence, $N(f,p^2)=p^2-p+1$.

When $r=3$, we need to consider the third term in the Taylor expansion $(2.1)$. That is,
$$f(a+lp)\equiv f(a)+f'(a)lp+\frac{f^{(2)}(a)}{2}(lp)^2(\bmod\ p^3).$$
Since $\frac{f^{(2)}(X)}{2}=(2p-1)(p-1)X^{2p-3}$ and $2p-3\geq 1$,
$\frac{f^{(2)}(a)}{2}=(2p-1)(p-1)a^{2p-3}\equiv0(\bmod\ p)$ for any $a\equiv0(\bmod\ p)$. Then, for $b=a+lp$ with $a\equiv0(\bmod\ p)$ and $l\not\equiv0(\bmod\ p)$, since 
$\text{ord}_p(f'(a)lp)=2$,
$$f(b)\equiv f(a)+f'(a)lp+\frac{f^{(2)}(a)}{2}(lp)^2\equiv f(a)+f'(a)lp\not\equiv f(a)(\bmod\ p^3).$$
Hence, $(f\bmod\ p^3)$ is injective restricted on $R_3\setminus S_3(0)$ and 
$$N(f,p^3)=p^3-p^2+p.$$ 

The proof is finished.  $\Box$

{\bf Remark. } Lemma \ref{lem:3.4} shows that given a prime $p$, we can construct a unified non-permutation polynomial $f(X)\in\ZZ[X]$ such that the value set of $(f\bmod\ p^r)$ are maximal in all non-permutation polynomial functions for all $r\geq 2$.

Putting the above things together, we obtain

{\bf The Proof of Theorem 1.1} When $r=1$, if $f$ is not a permutation over $\ZZ/p\ZZ$, $V(f,p)$ has $(p-1)$ elements at most. And, $M(p)=p-1$ since every mapping on $\ZZ/p\ZZ$ is a polynomial function.

When $r\geq 2$, the polynomial $f(X)=X^{2p-1}+pX$ in Lemma \ref{lem:3.4} shows that $$M(p^r)\geq N(f,p^r)= p^{r-2}(p^2-p+1).$$

On the other hand, if $f\in\ZZ[X]$ is not a permutation polynomial over $\ZZ/p^r\ZZ$ for $r\geq 2$, then $f$ is not a permutation over $\ZZ/p\ZZ$ or
$f'(a)\equiv0(\bmod\ p)$ for some $a\in\ZZ$. Then, by 
Lemma \ref{lem:3.1} and \ref{lem:3.2},
$$M(p^r)\leq\max\{p^r-p^{r-1},p^r-p^{r-1}+p^{r-2}\}=p^{r-2}(p^2-p+1).$$

We obtain that $M(p^r)=p^{r-2}(p^2-p+1)$. $\Box$

\section{The value $M(m)$ for $m=p_1^{r_1}\cdots p_k^{r_k}$}

For a positive integer $m$ with prime decomposition
$m=p_1^{r_1}\cdots p_k^{r_k}$, there is an isomorphism of rings
$$\Phi:\ZZ/m\ZZ\simeq(\ZZ/p_1^{r_1}\ZZ)\times\cdots\times(\ZZ/p_k^{r_k}\ZZ)$$
by Chinese remainder theorem.

One has the following two lemmas
\begin {lem}\label{lem:4.1}For an integer $m=p_1^{r_1}\cdots p_k^{r_k}$ and a polynomial $f\in\ZZ[X]$, 
$$V(f\bmod\ m)\simeq V(f\bmod\ p_1^{r_1})\times\cdots\times V(f\bmod\ p_k^{r_k}),$$
and 
$$N(f,m)=N(f, p_1^{r_1})\times\cdots\times N(f, p_k^{r_k}).$$
\end{lem}	
{\bf Proof. }Since $a\equiv b(\bmod\ m)$ implies $a\equiv b(\bmod\ p_i^{r_i})$, $1\leq i\leq k$, the morphism $\Phi$ restricted on $V(f\bmod\ m)$ gives us a mapping
\begin{align*}\Phi:V(f\mod\ m)&\rightarrow V(f\bmod\ p_1^{r_1})\times\cdots\times V(f\bmod\ p_k^{r_k})\\ (f(a)\bmod m)&\mapsto(f(a)\bmod\ p_1^{r_1},\cdots,f(a)\bmod\ p_k^{r_k}),\end{align*}
which is injective by Chinese remainder theorem. 

On the other hand, for any $(b_1,\cdots,b_k)\in V(f\bmod\ p_1^{r_1})\times\cdots\times V(f\bmod\ p_k^{r_k})$, there is $(a_1,\cdots,a_k)\in\ZZ^k$ such that 
$$f(a_i)\equiv b_i(\bmod\ p_i^{r_i}), \text{ for all }1\leq i\leq k.$$
Also by Chinese remainder theorem, there exists an integer 
$a$ satisfying $a\equiv a_i(\bmod\ p_i^{r_i})$ for all $1\leq i\leq k$. Then, $f(a)\equiv f(a_i)(\bmod\ p_i^{r_i})$
 and $$(f(a)\bmod p_1^{r_1},\cdots,f(a)\bmod p_k^{r_k})=(b_1,\cdots,b_k).$$
Hence, $\Phi$ is also surjective.

We have that $\Phi$ is a bijection on $V(f\bmod\ m)$, the size of which is the product of the sizes of $V(f\bmod\ p_i^{r_i}), 1\leq i\leq k$.   $\Box$

In particular, $f$ is a permutation polynomial over $\ZZ/m\ZZ$ if and only if $f$ is a permutation polynomial over $\ZZ/p_i^{r_i}\ZZ$ for any $1\leq i\leq k$.

\

Meanwhile,
\begin {lem}\label{lem:4.2}For an integer $m=p_1^{r_1}\cdots p_k^{r_k}$ and integral polynomials $f_1,\cdots,f_k\in\ZZ[X], 1\leq i\leq k$, there exists a polynomial $f\in\ZZ[X]$ such that 
$$V(f\bmod\ m)\simeq V(f_1\bmod\ p_1^{r_1})\times\cdots\times V(f_k\bmod\ p_k^{r_k}).$$
\end{lem}	
{\bf Proof. }Applying Chinese remainder theorem on the coefficients of $f_i$, we can obtain a polynomial $f\in\ZZ[X]$ such that $f\equiv f_i(\bmod\ p_i^{r_i})$ for all $1\leq i\leq k$. Then, $$V(f\bmod\ p_i^{r_i})=V(f_i\bmod p_i^{r_i}), 1\leq i\leq k.$$  
By Lemma \ref{lem:4.1}, we have the result.  $\Box$

\

Furthermore, for a prime factor $p$ of $m$, set the number
 $$m(p)=\begin{cases}m\left(1-\frac{1}{p}\right), &\text{ord}_p(m)=1\\m\left(1-\frac{1}{p}+\frac{1}{p^2}\right) &\text{ord}_p(m)\geq 2
 \end{cases}.$$
Then,
\begin {thm}\label{thm:4.2}For an integer $m=p_1^{r_1}\cdots p_k^{r_k}$, let $p$ be the maximal prime factor of $m$. One has that
$$M(m)=m(p),$$
unless $m=2^r3$ for $r\geq 2$ and $M(m)=2^{r-2}3^2$ in this case.
\end{thm}	
{\bf Proof.} If $f$ is not a permutation polynomial over $\ZZ/m\ZZ$, there is some $1\leq i\leq k$ such that 
$f$ is not a permutation polynomial over $\ZZ/p_i^{r_i}\ZZ$ by Lemma \ref{lem:4.1}.
Using Theorem \ref{thm:1.1} and Lemma \ref{lem:4.2}, 
it is obvious that
$$M(m)=\max\{m(p_i)\ |\ 1\leq i\leq k\}.$$

Observe that both $1-\frac{1}{p}$ and $1-\frac{1}{p}+\frac{1}{p^2}$ are strictly ascending when $p$ is increasing. On the other hand, if $p'$ is another prime less than $p$, then $1-\frac{1}{p'}+\frac{1}{p'^2}>1-\frac{1}{p}$ holds only if $p=3,p'=2$. We obtain the last assertion. $\Box$

\end{document}